\documentclass[12pt, a4paper]{amsart}
\usepackage[top=1.5in, bottom=1.5in, left=1in, right=1in]{geometry}
\usepackage[utf8]{inputenc}

\usepackage{mathrsfs}
\usepackage{latexsym}
\usepackage{graphics}
\usepackage{graphicx}
\usepackage{hyperref, amsmath, amsthm, amsfonts, amscd, flafter,epsf}
\usepackage{amsmath,amsfonts,amsthm,amssymb,amscd,amsxtra}
\usepackage{amsthm}
\input amssym.def
\input amssym.tex

\usepackage{hyperref}
\usepackage[alphabetic]{amsrefs}

\newtheorem{theorem}{Theorem}

\newtheorem{corollary}{Corollary}
\newtheorem{lemma}{Lemma}

\theoremstyle{remark}

\theoremstyle{theorem}

\newcommand{\ud}{\mathrm{d}}

\newcommand{\abs}[1]{\left|#1\right|}

\newcommand{\lcm}{\mathrm{lcm}}

\renewcommand{\mod}[1]{{\ifmmode\text{\rm\ (mod~$#1$)}\else\discretionary{}{}{\hbox{ }}\rm(mod~$#1$)\fi}}

\newcommand{\chibar}{\overline \chi}
\newcommand{\qr}[2]{\genfrac{(}{)}{}{}{#1}{#2}}

\title{Averages of character sums}
\author{Jonathan W. Bober}
\address{Heilbronn Institute for Mathematical Research \\ Department of Mathematics, University of Bristol,
Bristol, United Kingdom}
\email{{\tt j.bober@bristol.ac.uk}}
\date{September 5, 2014}
\begin{document}

\maketitle
\begin{abstract}
We show that a short truncation of the Fourier expansion for a character sum
gives a good approximation for the average value of that character sum over
an interval.

We give a few applications of this result. One is that for any $b$ there
are infinitely many characters for which the sum up to $\approx aq/b$ is
$\gg q^{1/2} \log \log q$ for all $a$ relatively prime to $b$;
another is that if the least quadratic nonresidue
modulo $q \equiv 3 \pmod 4$ is large, then the character sum gets as large as
$(\sqrt{q}/\pi) (L(1, \chi) + \log 2 - \epsilon)$, and if $B$ is this
nonresidue, then there is a sum of length $q/B$ which has size
$(\sqrt{q}/\pi) (\log 2 - \epsilon)$.
\end{abstract}

\thispagestyle{empty}

\section{Introduction}

For a primitive Dirichlet character $\chi$ modulo $q$, the character sum
\[
    S_\chi(x) = \sum_{n \le x} \chi(n)
\]
is a periodic function with period $q$. The Fourier expansion of this sum,
\[
    S_\chi(qt) \doteq \frac{\tau(\chi)}{2\pi i}\sum_{\abs{n} \ne 0}
        \frac{\chibar(n)}{n}\big(1 - e(nt)\big)
\]
(exactly valid for any $t$ which is not a discontinuity of the function), is a
useful tool for its study. This sum can be shortened considerably without
losing much accuracy; an example quantitative version is
\[
    S_\chi(qt) = \frac{\tau(\chi)}{2\pi i}\sum_{0 < |n| < q}
        \frac{\chibar(n)}{n}\big(1 - e(nt)\big) + O(\log q).
\]
By trivially bounding the terms in this sum, for instance, we can arrive
at the well-known Pólya--Vinogradov inequality
\[
    S_\chi(x) \ll q^{1/2} \log q.
\]
If we think that for large enough $n$ the values of the character sum are
fairly randomly distributed, then we might believe that a much shorter
truncation of this sum is still an accurate approximation to the character
sum. We can construct characters for which the values of $\chi(p)$ can be
controlled for all primes $p \ll \log q$, but beyond this point it becomes
much more difficult to understand the values of the character, and it seems
natural to expect the behavior to become more and more ``random'', from which
we might expect something like
\[
    S_\chi(qt) \approx \frac{\tau(\chi)}{2\pi i}\sum_{0 < |n| < (\log q)^A}
        \frac{\chibar(n)}{n}\big(1 - e(nt)\big)
\]
for some $A$. On the other hand, we cannot truncate the sum too early, as there
are examples of characters $\chi$ and real numbers $t$ such that
\[
    \sum_{|n| \ge \log q} \frac{\chibar(n)}{n}\big(1 - e(nt)\big) \gg \log\log q.
\]
This difficulty in truncating the sum early is an obstruction in attempts
to improve on the Pólya--Vinogradov inequality, and currently it is only
through the assumption of the Generalized Riemann Hypothesis that any
completely general improvement on this bound is known (\cite{MV}).

In this paper we show that if we are interested in the average value of the
character sum in a small interval, instead of the value of the sum at an exact
point, then a short truncation of the Fourier expansion can in fact be quite
accurate. We prove the following theorem.
\begin{theorem}\label{thm-average}
    For a primitive character $\chi$ mod $q$,
    \[
        \frac{B}{2}\int_{\alpha - 1/B}^{\alpha + 1/B} S_\chi(tq) \ud t =
        A(\chi) + \frac{\tau(\chi)}{i\pi}\sum_{n < B} \frac{\chibar(n)}{n}
        f(2 \pi n \alpha) + O\big(q^{1/2}\big),
    \]
    where $f(x) = -\cos(x)$ if $\chi(-1) = -1$ and $f(x) = i\sin(x)$ otherwise,
    and
    \[
        A(\chi) = \frac{\big(1 - \chi(-1)\big)\tau(\chi)}{2\pi i} L(1, \chibar).
    \]
\end{theorem}
One immediate application of this is that for even characters an improvement of
the Pólya--Vinogradov inequality for sums of length $o(q)$ would yield an
improvement for sums of length $O(q)$ as well.
\begin{corollary}
    For a primitive character $\chi$ with $\chi(-1) = 1$,
    \[
        M(\chi) \ll q^{1/2} \log\log q
        + \max_{N \le q} \max_{M \le q/\log q}
        \abs{\sum_{N \le n < N + M} \chi(n)}.
    \]
\end{corollary}
In fact, we can get this corollary with a small enough constant so that when
it is combined with lower bounds for the size of $M(\chi)$, we obtain
character sums of length $< q/\log q$ which have size
$\gg q^{1/2} \log \log q$. Granville and Soundararajan give stronger results of
this form in \cite{large-character-sums}, but the same result does not
explicitly appear for characters of larger order.

\begin{corollary}\label{corollary-large-sums}
    For any $A < \frac{e^\gamma}{\sqrt{3}}$ there exist infinitely many
    primitive even $\chi$ mod $q$ of arbitrarily large order such that
\[
    \abs{\sum_{q/3 \le n \le q/3 + N} \chi(n)} \ge
        \big(1 + o(1)\big)\left(\frac{e^\gamma}{\pi \sqrt 3} - \frac{A}{\pi}\right)
        q^{1/2}\log\log q
\]
    for some $|N| \le q/(\log q)^A$.
\end{corollary}

We will next examine Paley's construction \cite{paley} of quadratic characters
with $M(\chi) \gg q^{1/2} \log\log q$. Paley's construction suggests that the
character sums he considers will get large when they have length $q/4$, but
his proof does not actually give this result. Using Theorem \ref{thm-average},
we will fall only slightly short of this, showing that with Paley's
construction the character sum gets large at $(1/4 + t)q$ for some
$t \ll 1/\log q$. Moreover, we can get such a result for any rational number.
\begin{theorem}\label{thm-large-sums}
    For any $b > 0$, fix any primitive $\psi \mod b$ of order $g$. There exist
    infinitely many $\chi \mod q$ with parity opposite to $\psi$ and order
    $\lcm(2,g)$ such that for every $a$ coprime to $b$
    \[
        \abs{S_\chi(\alpha q) - A(\chi)} \ge
         \frac{\sqrt{q}}{\pi \sqrt{b}}\big(1 + o(1)\big) \log\log q
    \]
    for some $\alpha$ satisfying $\abs{\alpha - a/b} \ll 1/\log q$.
    Moreover, if $\chi$ is odd and $b > 1$, then assuming the Generalized
    Riemann Hypothesis, we may take $A(\chi) \ll q^{1/2}$ for these $\chi$.
\end{theorem}
    When $\chi$ is even $A(\chi) = 0$, so $S_\chi(\alpha q)$ is large. When
    $\chi$ is odd, we will be using GRH to show that $L(1, \chi) \ll 1$
    for the characters we construct if $b = 1$. Without the GRH assumption, all
    we can conclude is that either $S_\chi(\alpha q)$ is particularly large or
    $L(1, \chi)$ is particularly large. When $b = 1$, we expect that $L(1,
    \chi)$ will be large.

    To prove this theorem we will examine what happens when $\chi(n) = \psi(n)$
    for all small $n$. As such, we expect that this is a prototypical example
    of what happens when (in the language of Granville and Soundararajan
    \cite{pretentious-characters}) a character of large conductor is
    pretentious to a character of small conductor and opposite parity.
    
    This theorem includes the existence of characters of arbitrary even
    order for which $M(\chi) \gg q^{1/2}\log\log q$, a result of
    Goldmakher and Lamzouri \cite{GLeven}. Our construction of these characters
    is the same as theirs, though our proof that the sum is large differs.
    Both results suffer from the defect that the maximum of these
    character sums may decrease as the order increases.

    In fact, there exist characters of arbitrarily large (increasing) order
    such that $M(\chi) \ge (e^\gamma/\pi + o(1))q^{1/2}\log\log q$ (this
    follows from \cite[Theorem 3]{pretentious-characters}), and we expect
    such characters to be odd, and hence have even order. However, current
    results seem to still leave open the possibility that for any $\epsilon >
    0$ there exist (large) even $g$ such that all characters of order $g$
    satisfy $M(\chi) \le \epsilon q^{1/2}\log\log q$.  On the other hand, any
    analogue of Theorem \ref{thm-large-sums} which allowed for odd order
    characters would need to have the order of $\chi$ tending to infinity (as
    long as we accept the Generalized Riemann Hypothesis; see \cite[Theorem
    4]{pretentious-characters} and \cite[Theorem 2]{mimicry}).

    It is possible to obtain a slightly better constant in the above theorem by
    employing a version of Vinogradov's trick (as in \cite{pvlqnr}, for example).
    The theorem is probably true with the lower bound replaced by
    \[
        \frac{e^\gamma \sqrt{q}}{\pi \sqrt{b}}\big(1 + o(1)\big)\log\log q,
    \]
    and $\alpha = a/b$. In fact, if $\chi$ mimics $\psi$ sufficiently closely
    and $b > 1$, then it should be the case that
    \begin{equation}\label{eq-conjecture}
        S_\chi(aq/b) \sim \psi(a)\tau\big(\overline\psi\big) \frac{e^\gamma
            \tau(\chi)}{i\pi b}\log\log q.
    \end{equation}

The proof of the above theorem will reveal some partial information about the
direction in which the large value points, but because we get the result by
averaging the character sum over an interval, this information seems difficult
to extract from a complex-valued sum. For quadratic characters, we can take
advantage of the fact that the character sum is real and obtain some
approximation of equation \eqref{eq-conjecture}. As an example, we prove
\begin{theorem}\label{thm-large-quadratic}
    Fix any primitive odd quadratic character $\psi \bmod b$. There exist
    infinitely many even quadratic characters $\chi \bmod q$ such that
    for all $a$ relatively prime to $b$
    \[
        \sum_{n < \alpha q} \chi(n)
            = \psi(a) \frac{\sqrt{q}}{\pi\sqrt{b}}\big(1 + o(1)\big)\log \log q
    \]
    for some $\alpha$ satisfying $\abs{\alpha - a/b} \ll 1/\log q$.
\end{theorem}

Finally, we will examine the case of prime $q \equiv 3 \pmod 4$ with $\qr{n}{q}
= 1$ for many small $n$. In this case we expect that the value of the
$L$-function at $1$ --- and hence the constant term in the Fourier expansion
--- is large. This makes the character sum $S_\chi(qt)$ large for almost all $t
\asymp 1$, and the character sum tends to be maximized somewhat close to the
central point.  However, the central point is a local minimum of the character
sum, and computations suggest that the maximum actually tends to occur at or
near $(B-1)/2B$, where $B$ is the smallest quadratic nonresidue modulo $q$.

What we manage to prove here, using the same method of proof as Theorem
\ref{thm-average} but examining this specific case more carefully, is that
$S_\chi(qt)$ does get a bit larger than its value at the central
point for some $t$ within $1/B$ of $1/2$.
\begin{theorem}\label{thm-lqnr}
    For any sequence of odd quadratic characters $\chi \bmod q$ such that
    the least $B$ with $\chi(B) = -1$ tends to infinity, there exists a
    $t_\chi \in [1/2 - 1/B, 1/2]$ for each $\chi$ such that
    \[
        S_\chi(qt_\chi)- S_\chi(q/2) >
            \frac{\log 2}{\pi}\big(1 + o(1)\big)q^{1/2}.
    \]
\end{theorem}

As with Theorem \ref{thm-large-sums}, the constant here is probably not the
best possible. The $(\log 2)/\pi$ comes from computing the average of the
character sum over an interval, and at some points we know that the
character sum is significantly smaller than the average value. So it is natural
to expect that the character sum also gets significantly larger than its
average value; perhaps the theorem would still be true if the lower
bound were doubled.

If we consider the extreme case where $B > \log q$, then it seems likely
that $S_\chi(q/2) \ge e^\gamma/\pi(1 + o(1)) \log\log q$, and so
$(\log 2)\sqrt{q}/\pi = o(S_\chi(q/2))$. However, even this small difference
may have a large effect on the distribution of $M(\chi)$.
(See \cite{distribution-of-max2}.)

We also note that Theorem \ref{thm-lqnr} has the following consequence: If
$\sum_{N \le n < N + M} \chi(n) = o(q^{1/2})$ for all $M < q^{1-\epsilon}$,
then the least quadratic nonresidue for $q \equiv 3 \pmod 4$ is
$\ll q^\epsilon$.

\section{Averages of character sums}

We begin by recalling the Fourier expansion of $S_\chi(tq)$. For our purposes
in this paper, it will be convenient make explicit the difference between odd
and even characters, so we will work with the form
\[
    \widetilde S_\chi(tq) = A(\chi)
        + \frac{\tau(\chi)}{i\pi}F(t, \chi),
\]
where
\[
    F(t, \chi) = -\sum_{n=1}^\infty \frac{\chibar(n)}{n} \cos(2 \pi n t)
\]
if $\chi$ is odd ($\chi(-1) = -1$) and
\[
    F(t, \chi) = i\sum_{n=1}^\infty \frac{\chibar(n)}{n} \sin(2 \pi n t)
\]
if $\chi$ is even ($\chi(1) = 1$). As the character sum is not a continuous
function, the Fourier series does not always converge to the value of the
character sum, so here we are defining
\[
    \widetilde S_\chi(x)
        = \lim_{\epsilon \rightarrow 0}
            \frac{1}{2} \Big(S_\chi(x + \epsilon) + S_\chi(x - \epsilon)\Big).
\]
However, if we integrate the character sum, these discontinuities are
unimportant, and hence
\[
    \int_a^b S_\chi(tq) \ud t= (b-a)A(\chi)
                            + \frac{\tau(\chi)}{\pi}\int_a^b F(t, \chi) \ud t.
\]
We'll be particularly interested in the case where we integrate to compute
the average value over an interval a little bit shorter than $q$. We examine
\[
    \frac{B}{2}\int_{\alpha - 1/B}^{\alpha + 1/B} S_\chi(tq) \ud t.
\]
For this average, it will turn out that we can truncate the Fourier expansion
at $B$ and still obtain a good approximation. We recall Theorem
\ref{thm-average}, stated in the introduction.
\newtheorem*{thm:average}{Theorem \ref{thm-average}}
\begin{thm:average}
For a primitive character $\chi$ mod $q$,
\[
    \frac{B}{2}\int_{\alpha - 1/B}^{\alpha + 1/B} S_\chi(tq) \ud t = A(\chi) + \frac{\tau(\chi)}{i\pi}\sum_{n < B} \frac{\chibar(n)}{n} f(2 \pi n \alpha) + O\big(q^{1/2}\big),
\]
where $f(x) = -\cos(x)$ if $\chi(-1) = -1$ and $f(x) = i\sin(x)$ otherwise, and
\[
    A(\chi) = \frac{\big(1 - \chi(-1)\big)\tau(\chi)}{2\pi i} L(1, \chibar).
\]
\end{thm:average}
\begin{proof}
We wish to compute
\begin{equation} \label{eq:average1}
    \frac{B}{2}\int_{\alpha - 1/B}^{\alpha + 1/B} S_\chi(tq) \ud t
        = A(\chi) + \frac{B}{2}\frac{\tau(\chi)}{\pi}
            \int_{\alpha - 1/B}^{\alpha + 1/B} F(t, \chi) \ud t.
\end{equation}
We can evaluate the integral of $F(t, \chi)$ directly, integrating
term-by-term. For brevity we stick to the case where $\chi$ is odd (when
$\chi$ is even the proof is nearly identical, swapping sines and cosines).
We have
\[
\begin{split}
 \frac{B}{2}\int_{\alpha - 1/B}^{\alpha + 1/B} F(t, \chi) \ud t
  &= -\frac{B}{2}\sum_{n=1}^\infty \frac{\chibar(n)}{n} \int_{-1/B}^{1/B} \cos(2 \pi n t - 2\pi n \alpha) \ud t \\
  &= -\frac{B}{2}\sum_{n=1}^\infty \frac{\chibar(n)}{n} \int_{-1/B}^{1/B} \sin(2 \pi n t) \sin(2 \pi n \alpha)
        + \cos(2 \pi n t) \cos(2 \pi n \alpha) \ud t \\
  &= -B\sum_{n=1}^\infty \frac{\chibar(n)}{n} \cos(2 \pi n \alpha) \frac{\sin(2 \pi n/B)}{2\pi n}.
\end{split}
\]
The advantage of this process is that the sum now converges absolutely.
The tail of the sum is bounded in size by
\[
    \frac{B}{2\pi} \sum_{n=N+1}^\infty \frac{1}{n^2} < \frac{B}{2\pi N},
\]
so we truncate the sum at $B$ to obtain
\[
    \frac{B}{2}\int_{\alpha - 1/B}^{\alpha + 1/B} F(t, \chi) \ud t =
        -\sum_{n < B} \frac{\chibar(n)}{n} \cos(2 \pi n \alpha) \frac{\sin(2 \pi n/B)}{2\pi n/B} + O(1).
\]
Finally, using the approximation $\sin(x)/x = 1 + O(x^2)$, we have
\[
\begin{split}
    \frac{B}{2}\int_{\alpha - 1/B}^{\alpha + 1/B} F(t, \chi) \ud t
        &=-\sum_{n < B} \frac{\chibar(n)}{n} \cos(2 \pi n \alpha)
            + O\left(\frac{1}{B^2}\sum_{n < B}n\right) + O(1) \\
        &=-\sum_{n < B} \frac{\chibar(n)}{n} \cos(2 \pi n \alpha) + O(1).
\end{split}
\]
Putting this back into Equation \eqref{eq:average1} yields the theorem.
\end{proof}

\section{Large character sums at any point}
The existence of characters $\chi$ such that $M(\chi) \gg q^{1/2}\log\log q$
was first discovered by Paley (\cite{paley}). Paley constructed such characters
explicitly by using quadratic reciprocity to construct arbitrarily large
$q \equiv 1 \pmod 4$ such that
\[
    \qr{n}{q} = \left\{
        \begin{array}{rl}
            1  & \textrm{if $n \equiv 1 \pmod 4$} \\
            -1 & \textrm{if $n \equiv 3 \pmod 4$}
        \end{array}
        \right.
\]
for all $n \ll \frac{1}{4}\log q$. For $\chi = \qr{\cdot}{q}$, the Fourier
expansion for $S_\chi(q/4)$ is
\[
    S_\chi(q/4) = \frac{q^{1/2}}{\pi} \sum_{0 \le n < (\log q)/8}
     \frac{1}{2n + 1}
     + \frac{q^{1/2}}{\pi} \sum_{n \ge (\log q)/8} \frac{(-1)^n\chi(2n + 1)}{n}.
\]
Thus, we should probably expect $S_\chi(q/4)$ to be of size
$\gg q^{1/2}\log\log q$; however, the tail of this sum presents some
difficulties, and Paley proved his result in a manner that does not actually
give information about spot where $S_\chi(tq)$ is large. By applying Theorem
\ref{thm-average} we can average out the effects of the tail
of the sum to get the result that this sum is large near $t = 1/4$.
Moreover, by generalizing Paley's construction we can construct
even characters $\chi$ such that $S_\chi(tq)$ gets large near any rational
number $t$.

Paley's construction makes use of the fact that the sum
\[
    \sum_{n < x} \qr{-4}{n}\frac{\sin(2\pi n/4)}{n} \gg \log x.
\]
Here $\qr{-4}{n}$ is an odd character modulo $4$. More generally it turns out
that if $\chi \bmod q$ is an odd character, then
\[
    \sum_{n < x} \chi(n)\frac{\sin(2\pi n a/q)}{n} \gg \log x,
\]
while if $\chi$ is even,
\[
    \sum_{n < x} \chi(n)\frac{\cos(2\pi n a/q)}{n} \gg \log x.
\]
This is proved in the following lemma.
\begin{lemma}
For a primitive character $\chi$ mod $q$ and an $a$ relatively prime to $q$,
if $\chi(-1) = 1$, then
\[
    \sum_{n < x} \frac{\chi(n)}{n} \cos(2 \pi an/q)
        = \frac{\chibar(a)\tau(\chi)}{q}\log x + O(1),
\]
while if $\chi(-1) = -1$, then
\[
    \sum_{n < x} \frac{\chi(n)}{n} \sin(2 \pi an/q)
        = \frac{\chibar(a)\tau(\chi)}{iq}\log x + O(1).
\]
\end{lemma}
\begin{proof}
We may assume that $x$ is a multiple of $q$, since
$\sum_{n={mq}}^{(m+1)q} \frac{1}{n} \le 1/m$. We treat the case of $\chi(-1) = 1$;
the second case is basically identical.

We have
\[
\begin{split}
    \sum_{n < x} \frac{\chi(n)}{n} \cos(2 \pi an/q)
        &= \sum_{k=1}^{q-1} \sum_{n < x/q} \frac{\chi(nq + k)}{nq + k} \cos(2 \pi a(nq + k)/q) \\
        &= \sum_{k=1}^{q-1} \sum_{n < x/q} \frac{\chi(k)}{nq + k}\cos(2 \pi ak/q).
\end{split}
\]
Note that
\[
    \frac{1}{nq + k} = \frac{1}{nq} + \frac{k}{(nq)^2 + nqk},
\]
so the above sum is
\[
\sum_{n < x} \frac{\chi(n)}{n} \cos(2 \pi an/q) = \frac{A}{q} \sum_{n < x/q} \frac{1}{n} + O(1),
\]
where
\[
    A = \sum_{k = 1}^{q} \chi(k) \cos(2 \pi ak/q).
\]
Note that
\[
    \chibar(a)\tau(\chi) = \sum_{k=1}^q \chi(n)e(ak/q) = A + i \sum_{k = 1}^{q} \chi(k) \sin(2 \pi ak/q);
\]
however, since $\chi$ is even and sine is odd, the sum on the right hand side vanishes,
and we have
\[
\begin{split}
\sum_{n < x} \frac{\chi(n)}{n} \cos(2 \pi an/q)
    &= \frac{\chibar(a)\tau(\chi)}{q} \sum_{n < x/q} \frac{1}{n} + O(1) \\
    &= \frac{\chibar(a)\tau(\chi)}{q} \log x + O(1). \qedhere
\end{split}
\]

\end{proof}

Using this lemma and Theorem \ref{thm-average}, we can prove that there
are characters for which the character sum large gets near any rational number
that is not $0$ or $1/2$. We recall Theorem \ref{thm-large-sums}, which
was stated introduction.
\newtheorem*{thm-large-sums}{Theorem \ref{thm-large-sums}}
\begin{thm-large-sums}
    For any $b > 0$, fix any primitive $\psi \mod b$ of order $g$. There exist
    infinitely many $\chi \mod q$ with parity opposite to $\psi$ and order
    $\lcm(2,g)$ such that for every $a$ coprime to $b$
    \[
        \abs{S_\chi(\alpha q) - A(\chi)} \ge
         \frac{\sqrt{q}}{\pi \sqrt{b}}\big(1 + o(1)\big) \log\log q
    \]
    for some $\alpha$ satisfying $\abs{\alpha - a/b} \ll 1/\log q$.
    Moreover, if $\chi$ is odd and $b > 1$, then assuming the Generalized
    Riemann Hypothesis, we may take $A(\chi) \ll q^{1/2}$ for these $\chi$.
\end{thm-large-sums}
\begin{proof}
    A standard modification of Paley's construction gives the existence
    of a primitive odd quadratic character $\chi_1$ of arbitrarily large
    conductor $q$ such that $\chi_1(n) = 1$ for all $n < \frac{1}{2}\log q$.
    (See \cite[Proposition 2.1]{GLeven}, for example.) We
    can then take any primitive $\psi$ mod $b$ and consider the character $\chi
    = \chi_1 \psi$ mod $q$, where $q = b q_1$, which is primitive and has order
    $\lcm(2, g)$, where $g$ is the order of $\psi$. For this character we have
\[
    \frac{B}{2}\int_{a/b - 1/B}^{a/b + 1/B} S_\chi(qt) \ud t - A(\chi)
        = \frac{\tau(\chi)}{i \pi}\sum_{n < B} \frac{\chibar(n)}{n}f(2\pi n a/b) + O\big(q^{1/2}\big),
\]
where again $f(x) = -\cos(x)$ if $\chi(-1) = -1$ and $f(x) = i\sin(x)$ otherwise.
Choose $B = \frac{1}{2}\log q_1 = \log q - \log b$. In this range,
$\chibar_1(n) = 1$, so we have
\[
    \frac{B}{2}\int_{a/b - 1/B}^{a/b + 1/B} S_\chi(qt) \ud t - A(\chi)
        = \frac{\tau(\chi)}{i \pi}\sum_{n < B} \frac{\overline{\psi}(n)}{n}f(2\pi n a/b) + O\big(q^{1/2}\big).
\]
Then by the previous lemma we have
\[
    \frac{B}{2}\int_{a/b - 1/B}^{a/b + 1/B} S_\chi(qt) \ud t - A(\chi)
        = \psi(a) \frac{\tau(\chi)\tau\left(\overline{\psi}\right)}{ib \pi} \log B + O\big(q^{1/2}\big).
\]
We conclude that for some $\alpha \in [a/b - 1/B, a/b + 1/B]$, $S(\chi)$ must
get at least as large as this average value, which has size
\[
    \frac{\sqrt{q}}{\pi \sqrt{b}} \log\log q + O\big(q^{1/2}\big).
\]
For the extra statement about the size of $A(\chi)$ when $\chi$ is odd and
$b > 1$, we appeal to Littlewood's GRH-conditional result \cite{littlewood} that
\[
  L(1, \chi) \sim \prod_{p < (\log q)^2}\left(1 - \frac{\chi(p)}{p}\right)^{-1}.
\]
Note that
\[
    \prod_{\log q < p < (\log q)^2}\left(1 - \frac{1}{p}\right)^{-1} \sim 2,
\]
so
\[
  L(1, \chi) \ll 2 \prod_{p < (\log q)}\left(1 - \frac{\chi(p)}{p}\right)^{-1}
    = 2 \prod_{p < (\log q)}\left(1 - \frac{\psi(p)}{p}\right)^{-1}
    \sim 2 L(1, \psi),
\]
which is bounded as $q \rightarrow \infty$.
\end{proof}

\begin{proof}[Proof of Theorem \ref{thm-large-quadratic}]
In the above proof, if we take $\psi$ to be quadratic, then $\chi$ will
be quadratic as well, so we know the signs of the Gauss sums. For
an even character we obtain
\[
    \frac B 2 \int_{a/b - 1/B}^{a/b + 1/B} S_\chi(qt) \ud t
        = \psi(a) \frac{\sqrt{q}}{\pi \sqrt{b}} \log\log q  + O\big(q^{1/2}\big).
\]
In this case, since $S_\chi(qt)$ is real and close enough to a continuous
function, it must be within $1/2$ of its average value at some point in the
interval, which gives Theorem \ref{thm-large-quadratic}.
\end{proof}

\section{Large character sums over short intervals}
The error term in Theorem \ref{thm-average} is too large to be of much use when
considering sums with length smaller than $q^{1-\epsilon}$. However, there is a
range of lengths between $q$ and $q^{1 - \epsilon}$ where we can say something,
and by combining known lower bounds for character sums with stronger upper
bounds for averaged sums, we can assert the existence of sums of size
$\asymp q^{1/2}\log\log q$ over intervals of length $o(q)$, which is our
Corollary \ref{corollary-large-sums} in the introduction of this paper.

Some results of this type for real characters are obtained by Granville and
Soundararajan in \cite{large-character-sums} (see Theorem 11, for example) for
character sums over initial intervals. Granville and Soundararajan's results
are stronger than what we obtain from Theorem \ref{thm-average}, but our
proof is very simple and applies to any character for which the sum up to $q/3$
is large, so it is not restricted to real characters.

\begin{proof}[Proof of Corollary \ref{corollary-large-sums}]
For a primitive even character, Theorem \ref{thm-average} trivially gives
\[
    \begin{split}
        \abs{\frac{B}{2}\int_{1/3 - 1/B}^{1/3 + 1/B} S_\chi(tq)\ud t}
        &\le \frac{q^{1/2}}{\pi}\sum_{n < B} \frac{1}{n} + O\bigl(q^{1/2}\bigr) \\
        &\le \frac{q^{1/2}\log B}{\pi} + O\bigl(q^{1/2}\bigr).
    \end{split}
\]
On the other hand, there exist primitive even characters such that
\[
    \abs{S_\chi(q/3)} \ge \big(1 + o(1)\big)\frac{e^\gamma}{\pi \sqrt 3}
        q^{1/2}\log\log q.
\]
(See \cite[Equation 1.8]{pretentious-characters}.)
This means that for such characters, whenever
\[
    \log B < \frac{e^\gamma}{\sqrt 3} \log\log q,
\]
the average value of the character sum over this interval is significantly
less than the maximum of the character sum on this interval. So for such $B$,
there exist sums of length $\le q/B$ which are large. If we take
$B = (\log q)^A$, then we find that for $A < \frac{e^\gamma}{\sqrt 3}$, there
exists an $N$ with $|N| \le q/(\log q)^A$ such that
\[
    \abs{\sum_{q/3 \le n \le q/3 + N} \chi(n)} \ge
        (1 - o(1))\left(\frac{e^\gamma}{\pi \sqrt 3} - \frac{A}{\pi}\right)
        q^{1/2}\log\log q.
\]
\end{proof}

\section{A look at odd characters with large least nonresidue}

We now use the same methods as above to examine what happens with a quadratic
character modulo a prime $q \equiv 3 \pmod 4$ such that all small primes are
quadratic residues modulo $q$. For such characters, we expect that $L(1, \chi)$
is exceptionally large.  The central value of the character sum is
$\frac{\sqrt{q}}{\pi}L(1, \chi)$, so the character sum should be quite large at
this point; this central point, however, will actually be a local minimum of
the character sum. In fact, computations tend to reveal that the maximum of
this sum occurs at or near $(B-1)/2B$, where $B$ is the least quadratic
nonresidue modulo $q$.

Recall the Fourier expansion for an odd quadratic character $\chi \bmod q$
with $\chi(2) = 1$:
\[
    \widetilde S_\chi(tq) = \frac{\sqrt{q}}{\pi} L(1, \chi)
        + \frac{\sqrt{q}}{\pi} F(t, \chi)
\]
where
\[
    F(t, \chi) = \sum_{n=1}^\infty \frac{\chi(n)}{n} \cos(2\pi nt).
\]
If $\chi(n) = 1$ for all $n < B$, then for $t = 1/2$, the sum in the expansion
starts
\[
    -\sum_{n=1}^{B-1} \frac{(-1)^{n}}{n} \approx \log 2.
\]
However, for $t = 1/2$ we know that that
\[
    \sum_{n \le q/2} \chi(n) = \frac{\sqrt{q}}{\pi} L(1, \chi)
\]
exactly, which means that
\[
F(1/2, \chi) = 0,
\]
and thus the terms in the sum for $n \ge B$ completely cancel out the
initial contribution from the terms where $\chi(n) = 1$. However, as with
Theorem \ref{thm-average}, if we instead
consider the average value of this sum in an interval around $t = 1/2$, the
rapid oscillation of the cosines for large $n$ will cause these terms to have
a much smaller influence.

We consider the average of the character sum around $1/2$, computing the same
integral as in Theorem \ref{thm-average}.
\begin{align*}
\frac{B}{2}\int_{\frac{1}{2} - \frac{1}{B}}^{\frac{1}{2} + \frac{1}{B}}
    F(t, \chi) \ud t
  &= \frac{B}{2}\int_{-\frac{1}{B}}^{\frac{1}{B}}
        \sum_{n=1}^\infty \frac{\chi(n)}{n}(-1)^{n+1} \cos(2 \pi n t) \ud t \\
  &= B\sum_{n=1}^q \frac{(-1)^{n+1} \chi(n)}{n} \frac{\sin(2\pi n/B)}{2\pi n}.
\end{align*}
Instead of truncating now at $B$, we choose some (large) parameter $A$
and truncate the sum at $AB$, with an error of at most
$\frac{B}{2\pi} \cdot \frac{1}{AB}$, so that
\[
\frac{B}{2}\int_{\frac{1}{2} - \frac{1}{B}}^{\frac{1}{2} + \frac{1}{B}} F(t, \chi) \ud t
 = B\sum_{n=1}^{AB} \frac{(-1)^{n+1} \chi(n)}{n} \frac{\sin(2\pi n/B)}{2\pi n} + O(1/A).
\]
Now, since asymptotically $100\%$ of the numbers up to $AB$ are $(B-1)$-smooth as
$B \rightarrow \infty$, this sum will be well approximated by $\chi(n)$
replaced by $1$. Specifically, since all of these numbers are $\ge B$, each
one contributes at most $\frac{1}{B^2}$ to the sum and an easy count reveals
that the set of $n$ in this range such that $\chi(n) \ne 1$ is very small.
\begin{lemma}
The number of integers up to $AB$ that are not $(B-1)$-smooth is $\ll \frac{A^2
B}{\log AB}$.
\end{lemma}
\begin{proof}
Suppose that $B \le n \le AB$ is divisible by some prime $p \ge B$. Then $n/p
\le A$.  This means that any number in the range $B \le n \le AB$ which is not
$B$-smooth is the product a prime in that range and a number less than $A$.
There are $\ll AB/\log(AB)$ primes in this range, and $A$ integers $< A$,
which gives the bound.
\end{proof}

Using this lemma, we see that if we replace $\chi(n)$ by $1$ in the sum
up to $AB$, the error will be at most
\[
    B \cdot \frac{1}{B^2} \cdot \frac{A^2 B}{\log AB} = \frac{A^2}{\log AB}.
\]
After replacing $\chi(n)$ by $1$, we can also extend the sum to infinity
at an error of at most $O(1/A)$, so we have
\[
\frac{B}{2}\int_{\frac{1}{2} - \frac{1}{B}}^{\frac{1}{2} + \frac{1}{B}} F(t, \chi) \ud t
 = B\sum_{n=1}^{\infty} \frac{(-1)^{n+1}}{n} \frac{\sin(2\pi n/B)}{2\pi n}
 + O\left(\frac{A^2}{\log AB}\right)
 + O(1/A).
\]
Now as $B$ gets large we see that this sum approaches
\[
    \lim_{x \rightarrow 0} \sum_{n=1}^\infty \frac{(-1)^{n+1}}{n} \frac{\sin(nx)}{nx}
         = \sum_{n=1}^\infty \frac{(-1)^{n+1}}{n} = \log 2.
\]
This proves Theorem \ref{thm-lqnr}, which we recall here.
\newtheorem*{thm-lqnr}{Theorem \ref{thm-lqnr}}
\begin{thm-lqnr}
    For any sequence of odd quadratic characters $\chi \bmod q$ such that
    the least $B$ with $\chi(B) = -1$ tends to infinity, there exists a
    $t_\chi \in [1/2 - 1/B, 1/2]$ for each $\chi$ such that
    \[
        S_\chi(qt_\chi)- S_\chi(q/2) >
            \frac{\log 2}{\pi}\big(1 + o(1)\big)q^{1/2}.
    \]
\end{thm-lqnr}

\end{document}